\documentclass{article}

\usepackage[top=20truemm,bottom=20truemm,left=20truemm,right=20truemm]{geometry}
\usepackage{authblk}
\usepackage{abstract}
\usepackage{amsmath}
\usepackage{amssymb}
\usepackage{hyperref}
\usepackage[dvipdfmx]{graphicx}
\usepackage{color}
\usepackage{diagbox}

\author[1]{Takuya Tsuchiya\thanks{\href{mailto:tatsuchi@eco.meijigakuin.ac.jp}{%
      \nolinkurl{tatsuchi@eco.meijigakuin.ac.jp}}}}
\author[2]{Makoto Nakamura}
\affil[1]{%
  Faculty of Economics, Meiji Gakuin University, Japan
}
\affil[2]{%
  Department of Pure and Applied Mathematics,
  Graduate School of Information Science and Technology,
  The University of Osaka, Japan
}

\title{Quantitative evaluations of stability and convergence
  for solutions of semilinear Klein--Gordon equation}

\begin{document}

\maketitle

\begin{abstract}
  We perform some simulations of the semilinear Klein--Gordon equation with a
  power-law nonlinear term and propose each of the quantitative evaluation
  methods for the stability and convergence of numerical solutions.
  We also investigate each of the thresholds in the methods by varying the
  amplitude of the initial value and the mass, and propose appropriate values.
\end{abstract}

\section{Introduction}

Many natural phenomena are expressed by (nonlinear) hyperbolic equations.
We are strongly interested in the behavior of the asymptotic solutions of the
equations in time.
In addition, the properties of the solutions should be changed in a curved
spacetime since the differential operator is affected by the curvature of
spacetime (e.g.\cite{Wald}).
We adopt the Klein--Gordon equation as the hyperbolic equation since it can be
applied to a curved spacetime.
Some analytical results of the equation in the de Sitter spacetime, which is one
of the curved spacetimes, have been reported \cite{Yagdjian-Galstian,
  Nakamura-2014-JMAA}.
Regarding the numerical study of the equation, we have reported the numerical
solutions of the equation in the de Sitter spacetime using the
structure-preserving scheme \cite{Tsuchiya-Nakamura-2019-JCAM}, suggested some
discrete equations constructed using the structure-preserving scheme
\cite{Tsuchiya-Nakamura-2022-ISAAC}, and investigated the reasons for the
difference in stability between the discrete equations \cite{%
  Tsuchiya-Nakamura-2023-JSIAML}.

In (partial) differential equations, the stability and convergence of numerical
solutions are necessary for the correctness of the solutions.
Although we have proposed highly accurate numerical solutions for the
semilinear Klein--Gordon equation \cite{Tsuchiya-Nakamura-2019-JCAM,
  Tsuchiya-Nakamura-2022-ISAAC, Tsuchiya-Nakamura-2023-JSIAML}, we have not
quantitatively evaluated the stability and convergence of the solutions.
In this paper, we propose some quantitative evaluation methods for the
stability and convergence of the solutions for the semilinear Klein--Gordon
equation in the flat spacetime.

Indices such as $(i, j, \dots)$ run from 1 to $n$, where $n$ is the spatial
dimension.
We use the Einstein convention of summation of repeated up--down indices in this
paper.

\section{Semilinear Klein--Gordon equation}

The semilinear Klein--Gordon equation with the power-law nonlinear term in the
flat spacetime is
\begin{align}
  -\dfrac{1}{c^2}\partial_t^2\phi + \delta^{ij}(\partial_i\partial_j\phi)
  - \dfrac{c^2m^2}{\hbar^2}\phi = \lambda |\phi|^{p-1}\phi,
  \label{eq:semi-KGeq}
\end{align}
where $\phi$ is the dynamical variable, $\delta^{ij}$ is the Kronecker delta,
$m$ is the mass, $c$ is the speed of light, $\hbar$ is the Dirac constant, $p$
is an integer larger than 2, and $\lambda$ is a constant and has a physical
dimension of $1/(\text{length})^2$.
When performing numerical calculations, the canonical form is preferable since
it is a system of first-order equations in time.
Moreover, it is easy to confirm the accuracy of the numerical calculations since
there is a constraint with respect to time such as the total Hamiltonian.
The Hamiltonian density of \eqref{eq:semi-KGeq} is given by
\begin{align}
  \mathcal{H}
  &=
  \dfrac{L_0}{2}\left(\dfrac{\psi^2}{L_0^2}
  + \delta^{ij}(\partial_i\phi)(\partial_j\phi)
  + \dfrac{c^2m^2\phi^2}{\hbar^2}
  + \dfrac{2\lambda|\phi|^{p+1}}{p+1}\right),
  \label{eq:Hamiltonian}
\end{align}
where $L_0$ is a constant value that makes the physical dimension of
$\mathcal{H}$ into an energy dimension and $\psi$ is the canonical momentum of
$\phi$.
Then, the canonical equations of \eqref{eq:semi-KGeq} are
\begin{align}
  \dfrac{1}{c}\partial_t\phi
  &= \dfrac{1}{L_0}\psi,
  \label{eq:phi}
  \\
  \dfrac{1}{c}\partial_t\psi
  &=
  L_0\delta^{ij}(\partial_i\partial_j\phi)
  - \dfrac{L_0c^2m^2}{\hbar^2}\phi
  - L_0\lambda |\phi|^{p-1}\phi.
  \label{eq:psi}
\end{align}

The discretized equations of \eqref{eq:Hamiltonian}, \eqref{eq:phi}, and
\eqref{eq:psi} can be respectively defined as
\begin{align}
  \mathcal{H}{}^{(\ell)}_{(\boldsymbol{k})}
  &:= \frac{L_0}{2}\biggl(
  \dfrac{(\psi^{(\ell)}_{(\boldsymbol{k})})^2}{L_0^2}
  + \delta^{ij}
  (\widehat{\delta}^{\langle1\rangle}_i
  \phi^{(\ell)}_{(\boldsymbol{k})})(\widehat{\delta}^{\langle1\rangle}_j
  \phi^{(\ell)}_{(\boldsymbol{k})})
  + \dfrac{c^2m^2}{\hbar^2}
  (\phi^{(\ell)}_{(\boldsymbol{k})})^2
  + \dfrac{2\lambda}{p+1}
  |\phi^{(\ell)}_{(\boldsymbol{k})}|^{p+1}
  \biggr),
  \label{eq:DiscreteHamiltonian}
\end{align}
\begin{align}
  &\frac{\phi^{(\ell+1)}_{(\boldsymbol{k})}
    - \phi^{(\ell)}_{(\boldsymbol{k})}}{c\Delta t}
  :=
  \dfrac{1}{2L_0}
  (\psi^{(\ell+1)}_{(\boldsymbol{k})} + \psi^{(\ell)}_{(\boldsymbol{k})}),
  \label{eq:DiscretePhi-I}
  \\
  &\frac{\psi^{(\ell+1)}_{(\boldsymbol{k})}
    - \psi^{(\ell)}_{(\boldsymbol{k})}}{c\Delta t}
  :=
  L_0\biggl(
  - \frac{\lambda}{p+1}
  \frac{|\phi^{(\ell+1)}_{(\boldsymbol{k})}|^{p+1}
    - |\phi^{(\ell)}_{(\boldsymbol{k})}|^{p+1}}{
    \phi^{(\ell+1)}_{(\boldsymbol{k})}- \phi^{(\ell)}_{(\boldsymbol{k})}}
  + \dfrac{\delta^{ij}
    \widehat{\delta}^{\langle1\rangle}_i
    \widehat{\delta}^{\langle1\rangle}_j
    (\phi^{(\ell+1)}_{(\boldsymbol{k})} + \phi^{(\ell)}_{(\boldsymbol{k})})
  }{2}
  - \frac{c^2m^2(\phi^{(\ell+1)}_{(\boldsymbol{k})}
    + \phi^{(\ell)}_{(\boldsymbol{k})})}{2\hbar^2}\biggr),
  \label{eq:DiscretePsi-I}
\end{align}
where ${}^{(\ell)}$ means the time index, ${}_{(\boldsymbol{k})}$ means the
space index, and $\boldsymbol{k}=(k_1,\dots,k_n)$.
$\widehat{\delta}^{\langle1\rangle}_i$ is the first-order central difference
operator defined as
\begin{align*}
  \widehat{\delta}^{\langle1\rangle}_{i}u^{(\ell)}_{(\boldsymbol{k})}
  := \dfrac{u^{(\ell)}_{(k_1,\dots,k_i+1,\dots,k_n)}
    - u^{(\ell)}_{(k_1,\dots,k_i-1,\dots,k_n)}}{2\Delta x^i}.
\end{align*}
Note that \eqref{eq:DiscreteHamiltonian}--\eqref{eq:DiscretePsi-I} are called
Form I in \cite{Tsuchiya-Nakamura-2023-JSIAML}.
Here, $\Delta x^i$ is the $i$ th grid range.
If $n=3$, for example, $\Delta x^1=\Delta x$, $\Delta x^2=\Delta y$, and
$\Delta x^3=\Delta z$.
The nonlinear term can be expressed as
\begin{align}
  \frac{|\phi^{(\ell+1)}_{(\boldsymbol{k})}|^{p+1}
    - |\phi^{(\ell)}_{(\boldsymbol{k})}|^{p+1}}{
    \phi^{(\ell+1)}_{(\boldsymbol{k})}-
    \phi^{(\ell)}_{(\boldsymbol{k})}}
  &= \{|\phi^{(\ell+1)}_{(\boldsymbol{k})}|^{p}
  + |\phi^{(\ell+1)}_{(\boldsymbol{k})}|^{p-1}|\phi^{(\ell)}_{(\boldsymbol{k})}|
  + \cdots
  + |\phi^{(\ell+1)}_{(\boldsymbol{k})}||\phi^{(\ell)}_{(\boldsymbol{k})}|^{p-1}
  + |\phi^{(\ell)}_{(\boldsymbol{k})}|^{p}\}
  \frac{|\phi^{(\ell+1)}_{(\boldsymbol{k})}|
    - |\phi^{(\ell)}_{(\boldsymbol{k})}|}{
    \phi^{(\ell+1)}_{(\boldsymbol{k})} - \phi^{(\ell)}_{(\boldsymbol{k})}}.
\end{align}
The total Hamiltonian $\int_{\mathbb{R}^n}\mathcal{H}\,dx^n$ at a discrete level
is preserved using \eqref{eq:DiscretePhi-I} and \eqref{eq:DiscretePsi-I}
\cite{Tsuchiya-Nakamura-2023-JSIAML}.

\section{Quantitative evaluations of stability and convergence}

In this paper, the word ``stable simulation'' means that no vibration occurs in
the waveform of $\phi$.
Moreover, to quantitatively evaluate stability, we define
\begin{align}
  d\phi^{(\ell)}_{(\boldsymbol{k})}:=\hat{s}^+_i
  \phi^{(\ell)}_{(\boldsymbol{k})}-\phi^{(\ell)}_{(\boldsymbol{k})}
\end{align}
and sum of $|d\phi^{(\ell)}_{(\boldsymbol{k})}|/(\text{grid})$ if
$d\phi^{(\ell)}_{(\boldsymbol{k})}$ satisfies the
condition
\begin{align}
  (\hat{s}^{+}_id\phi^{(\ell)}_{(\boldsymbol{k})}) d\phi^{(\ell)}_{%
    (\boldsymbol{k})}<0
  \label{eq:QuantitativeStab}
\end{align}
over $\boldsymbol{k}$, where $\hat{s}^{+}_{i}$ is the discrete operator that
shifts the space forward.
We call this value $SV_g$, which is determined for each grid, and consider the
simulation stable when $SV_g$ is less than or equal to the threshold
$\varepsilon_{\mathrm{s}}$.
$SV_g$ represents the number of vibrations per grid weighted by
$|d\phi^{(\ell)}_{(\boldsymbol{k})}|$.
We study the appropriate value of $\varepsilon_{\mathrm{s}}$ in Section
\ref{sec:Num}.

The word ``convergence'' means that $\phi$ approaches the exact solution with
an increasing number of grids.
To quantitatively determine convergence, we define the relative errors of
$\phi$:
\begin{align}
  CV_{g}(t)
  &:=\log_{10}\dfrac{\|\phi_{g}(x) - \phi_{\mathrm{G}}(x)\|_2}{%
    \|\phi_{\mathrm{G}}(x)\|_2},
\end{align}
where $\phi_g(x)$ is the value of $\phi$ for each grid number and
$\phi_{\mathrm{G}}(x)$ is that for the maximum grid number.
Since \eqref{eq:DiscretePhi-I} and \eqref{eq:DiscretePsi-I} have the
second-order convergence with respect to the number of grids \cite{%
  Tsuchiya-Nakamura-2023-JSIAML}, we define the difference in $CV_g(t)$ from
the second-order convergence as
\begin{align}
  DCV_{g}(t)
  &:=\left|CV_{\bar{\mathrm{G}}}(t) - CV_g(t)+
  \left(\log_2\dfrac{\bar{\mathrm{G}}}{g}\right)
  \left(\log_{10}4\right)\right|,
  \label{eq:QuantitativeConv}
\end{align}
where $\bar{\mathrm{G}}$ is the second largest grid number.
If $DCV_{g}(t)$ is less than or equal to the threshold
$\varepsilon_{\mathrm{c}}$, we decide that the convergence of the simulation is
satisfied.
We also study the appropriate value of $\varepsilon_{\mathrm{c}}$ in Section
\ref{sec:Num}.

\section{Numerical results
  \label{sec:Num}
}
In this section, we perform some simulations using the settings given below.
The initial conditions are set as $\phi(x,\bar{\boldsymbol{x}})=\phi(x)=
A\cos(2\pi x)$ and $\psi(x,\bar{\boldsymbol{x}})=\psi(x)=2\pi A\sin(2\pi x)$
for $x\in\mathbb{R}$ and $\bar{\boldsymbol{x}}\in\mathbb{R}^{n-1}$ with $A=2,3$
and $-1/2\leq x\leq 1/2$.
The boundary is periodic.
The physical parameters are $c=\hbar=L_0=1$.
The spatial dimension is $n=3$.
The grid ranges are $(\Delta x, \Delta t)=(1/250, 1/2500)$, $(1/500, 1/5000)$,
$(1/1000, 1/10000)$, $(1/2000, 1/20000)$, $(1/4000, 1/40000)$, and
$(1/8000, 1/80000)$.
The simulation time is $0\leq t\leq 1000$.
The number of exponents of the nonlinear term is $p=5$ and the coefficient
parameter is $\lambda=1$.
The mass $m$ ranges from $3.9$ to $4.2$ when $A=2$ and from $7.6$ to $8.2$ when
$A=3$.

\begin{figure}[htbp]
  \centering
  \begin{minipage}{0.32\hsize}
    \includegraphics[width=\hsize]{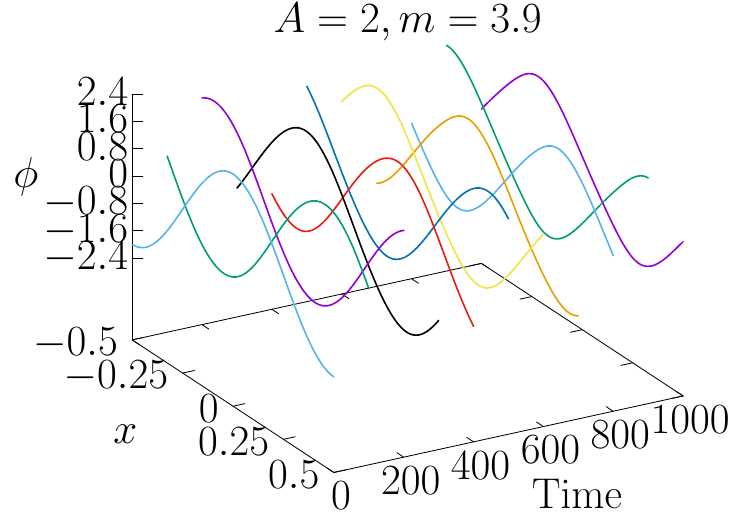}
  \end{minipage}
  \begin{minipage}{0.32\hsize}
    \includegraphics[width=\hsize]{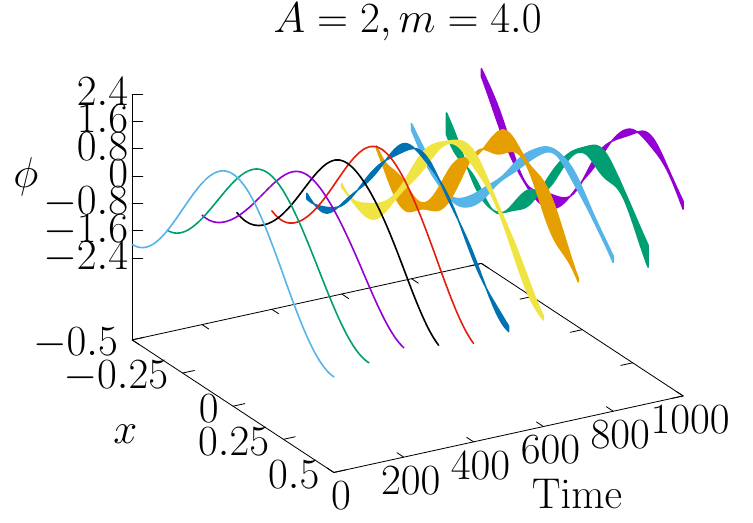}
  \end{minipage}
  \begin{minipage}{0.32\hsize}
    \includegraphics[width=\hsize]{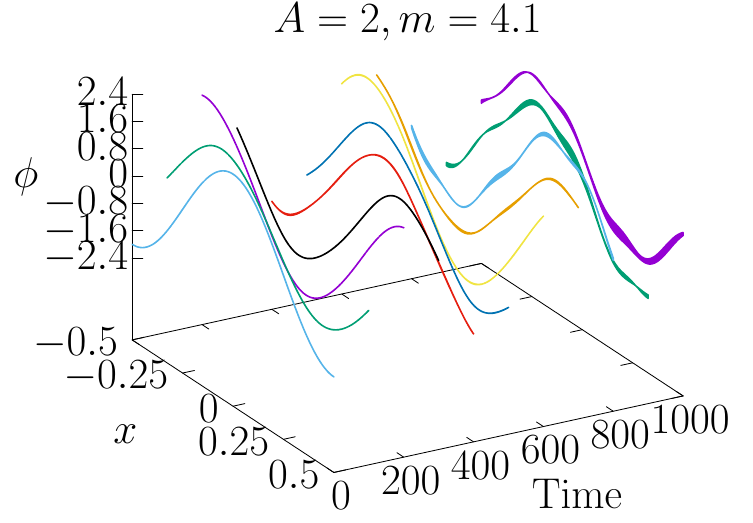}
  \end{minipage}
  \begin{minipage}{0.32\hsize}
    \includegraphics[width=\hsize]{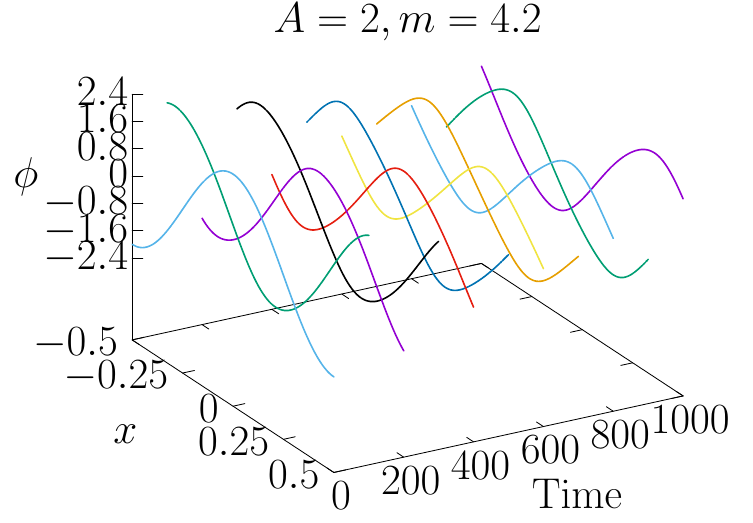}
  \end{minipage}
  \caption{
    $\phi$ with $A=2$, $m=3.9$ to $4.2$, and $8000$ grids.
    The top-left panel is for $m=3.9$, the top-center one is for $m=4.0$, the
    top-right one is for $m=4.1$, and the bottom one is for $m=4.2$.
    The vibration appears to occur at $t\geq 500$ for $m=4.0$ and at $t\geq 700$
    for $m=4.1$.
    \label{fig:WavesPhi8000_A2_m39_42}
  }
\end{figure}
Fig. \ref{fig:WavesPhi8000_A2_m39_42} shows $\phi$ with $A=2$ and $m=3.9$ to
$4.2$.
The vibration seems to occur at $t\geq 500$ for $m=4.0$ and at $t\geq 700$ for
$m=4.1$.
On the other hand, no vibration appears to occur for $m=3.9$ and $4.2$.
\begin{figure}[htbp]
  \centering
  \begin{minipage}{0.32\hsize}
    \includegraphics[width=\hsize]{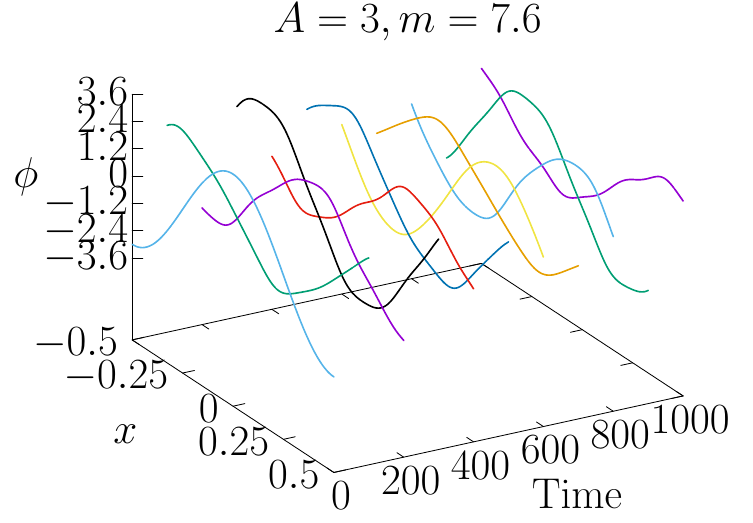}
  \end{minipage}
  \begin{minipage}{0.32\hsize}
    \includegraphics[width=\hsize]{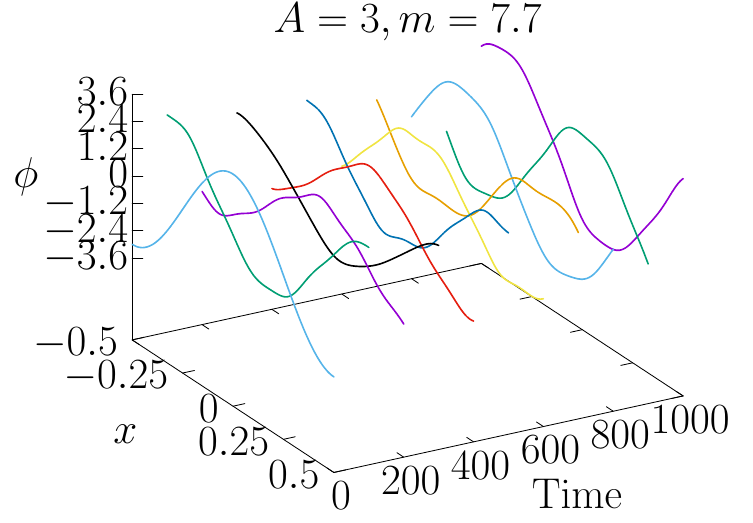}
  \end{minipage}
  \begin{minipage}{0.32\hsize}
    \includegraphics[width=\hsize]{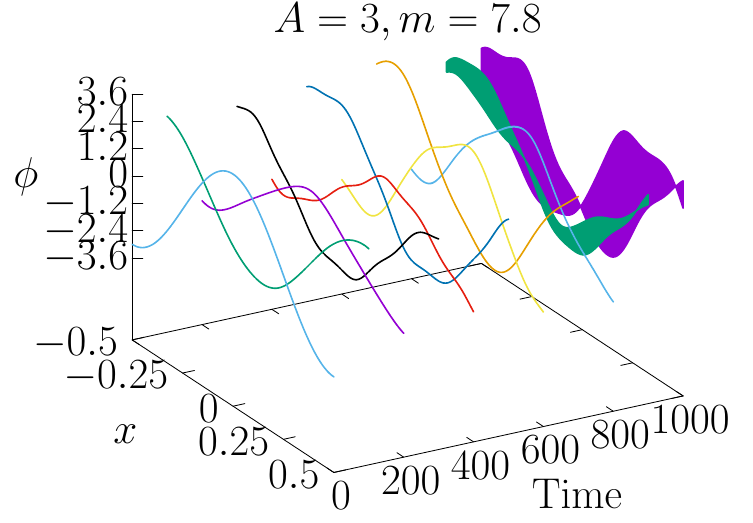}
  \end{minipage}
  \begin{minipage}{0.32\hsize}
    \includegraphics[width=\hsize]{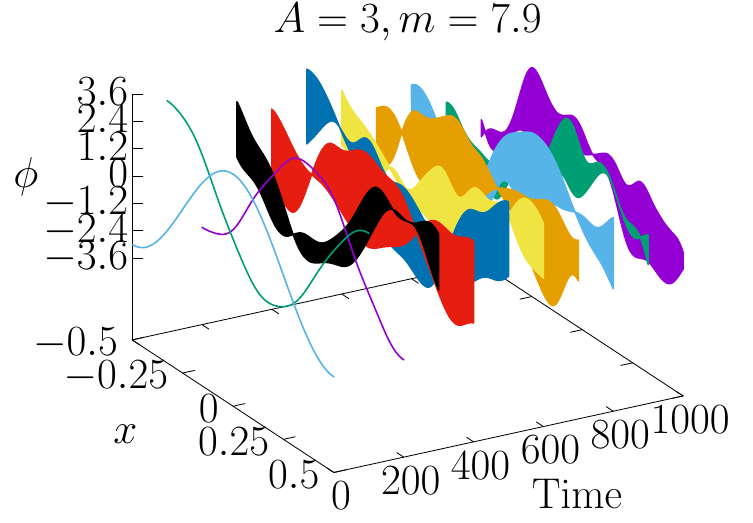}
  \end{minipage}
  \begin{minipage}{0.32\hsize}
    \includegraphics[width=\hsize]{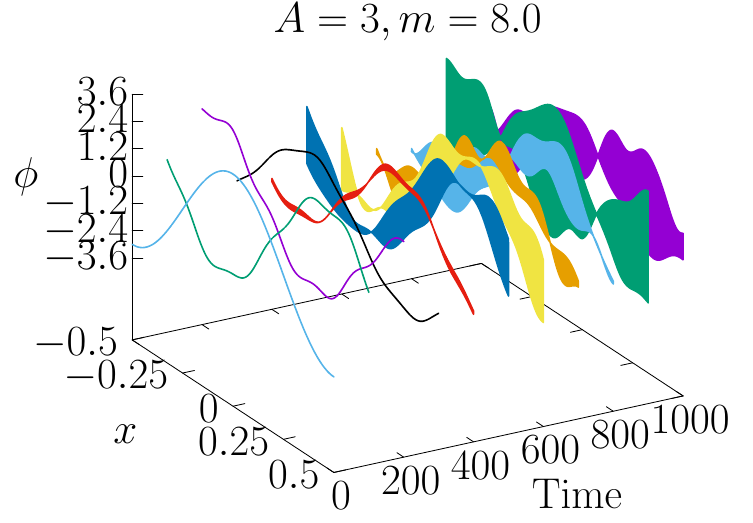}
  \end{minipage}
  \begin{minipage}{0.32\hsize}
    \includegraphics[width=\hsize]{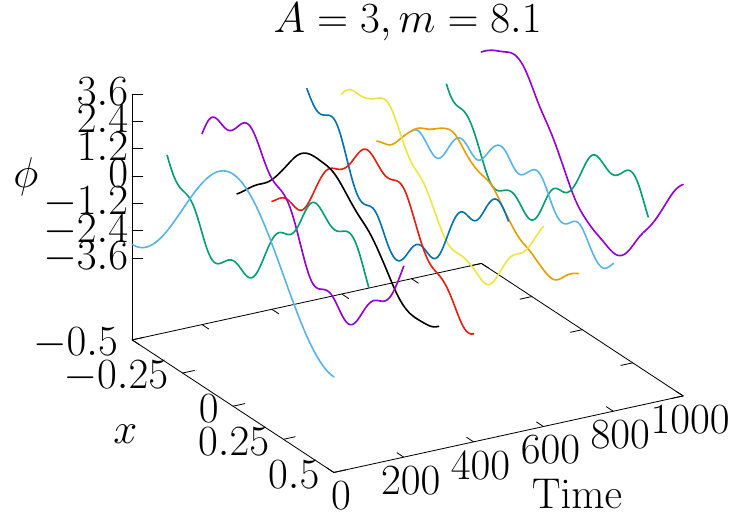}
  \end{minipage}
  \begin{minipage}{0.32\hsize}
    \includegraphics[width=\hsize]{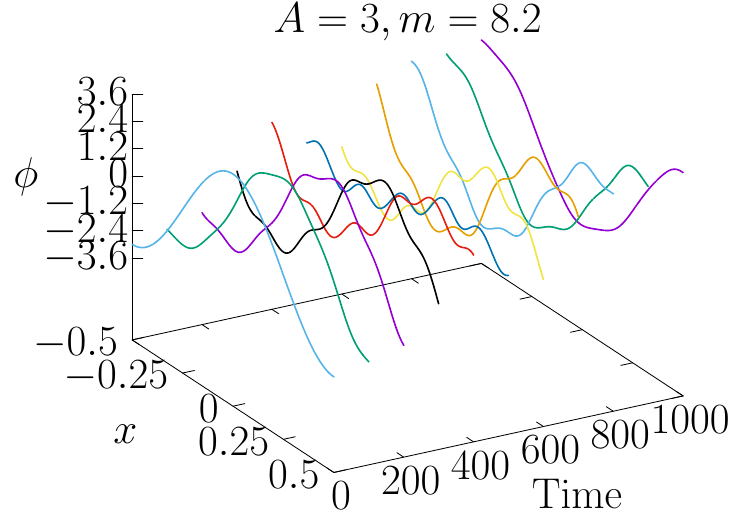}
  \end{minipage}
  \caption{
    $\phi$ with $A=3$, $m=7.6$ to $8.2$, and $8000$ grids.
    The top-left panel is for $m=7.6$, the top-center one is for $m=7.7$, the
    top-right one is for $m=7.8$, the center-left one is for $m=7.9$, the center
    one is for $m=8.0$, the center-right one is for $m=8.1$, and the bottom one
    is for $m=8.2$.
    The vibration appears to occur at $t\geq 900$ for $m=7.8$, at $t\geq 300$
    for $m=7.9$, and at $t\geq400$ for $m=8.0$.
    \label{fig:WavesPhi8000_A3_m76_82}
  }
\end{figure}
Fig. \ref{fig:WavesPhi8000_A3_m76_82} shows $\phi$ with $A=3$ and $m=7.6$ to
$8.2$.
The vibration seems to occur at $t\geq 900$ for $m=7.8$, at $t\geq 300$ for
$m=7.9$, and at $t\geq 400$ for $m=8.0$.
On the other hand, no vibrations appear to occur for $m=7.6$, $7.7$, $8.1$, and
$8.2$.

To quantitatively evaluate stability, we investigate the appropriate value
of $\varepsilon_{\mathrm{s}}$.
The value of $\varepsilon_{\mathrm{s}}$ is varied to $0.001$, $0.06$, $0.12$,
$0.18$, $0.24$, and $0.30$.
We summarize the results of the stability in Table \ref{tbl:stabA2A3}.
\begin{table}[htbp]
  \centering
  \begin{tabular}{|c||c|c|c|c|}
    \hline
    \multicolumn{5}{|c|}{$A=2$}\\
    \hline
    \diagbox[width=1.25cm]{$\varepsilon_{\mathrm{s}}$}{$m$}
    & $3.9$ & $4.0$ & $4.1$ & $4.2$ \\
    \hline
    $0.001$ & $\bigcirc$ & $417$ & $368$ & $\bigcirc$\\
    \hline
    $0.06$ & $\bigcirc$ & $474$ & $414$ & $\bigcirc$\\
    \hline
    $0.12$ & $\bigcirc$ & $488$ & $425$ & $\bigcirc$\\
    \hline
    $0.18$ & $\bigcirc$ & $496$ & $433$ & $\bigcirc$\\
    \hline
    $0.24$ & $\bigcirc$ & $497$ & $558$ & $\bigcirc$\\
    \hline
    $0.30$ & $\bigcirc$ & $504$ & $\bigcirc$ & $\bigcirc$\\
    \hline
  \end{tabular}
  \begin{tabular}{|c||c|c|c|c|c|c|c|}
    \hline
    \multicolumn{8}{|c|}{$A=3$}\\
    \hline
    \diagbox[width=1.25cm]{$\varepsilon_{\mathrm{s}}$}{$m$}
    & $7.6$ & $7.7$ & $7.8$ & $7.9$ & $8.0$ & $8.1$ & $8.2$\\
    \hline
    $0.001$ & $\bigcirc$ & $930$ & $725$ & $242$ & $343$ & $\bigcirc$
    & $\bigcirc$\\
    \hline
    $0.06$ & $\bigcirc$ & $\bigcirc$ & $845$ & $272$ & $393$ & $\bigcirc$
    & $\bigcirc$\\
    \hline
    $0.12$ & $\bigcirc$ & $\bigcirc$ & $887$ & $275$ & $400$ & $\bigcirc$
    & $\bigcirc$\\
    \hline
    $0.18$ & $\bigcirc$ & $\bigcirc$ & $890$ & $279$ & $404$ & $\bigcirc$
    & $\bigcirc$\\
    \hline
    $0.24$ & $\bigcirc$ & $\bigcirc$ & $891$ & $280$ & $404$ & $\bigcirc$
    & $\bigcirc$\\
    \hline
    $0.30$ & $\bigcirc$ & $\bigcirc$ & $892$ & $280$ & $406$ & $\bigcirc$
    & $\bigcirc$\\
    \hline
  \end{tabular}
  \caption{
    Time when $SV_{8000}>\varepsilon_{\mathrm{s}}$ for $A=2$ and $m=3.9$ to
    $4.2$, and for $A=3$ and $m=7.6$ to $8.2$.
    The left table is for $A=2$ and the right one is for $A=3$.
    $\bigcirc$ means that $SV_{8000}\leq\varepsilon_{\mathrm{s}}$ is always
    satisfied for $0\leq t\leq 1000$.
    The values represent the times when $SV_{8000}>\varepsilon_{\mathrm{s}}$.
    \label{tbl:stabA2A3}
  }
\end{table}
$\bigcirc$ means that $SV_{8000}\leq\varepsilon_{\mathrm{s}}$ is always
satisfied at $0\leq t\leq 1000$.
On the other hand, the values represent the times when $SV_{8000}>
\varepsilon_{\mathrm{s}}$.

\begin{figure}[htbp]
  \centering
  \begin{minipage}{0.32\hsize}
    \includegraphics[width=\hsize]{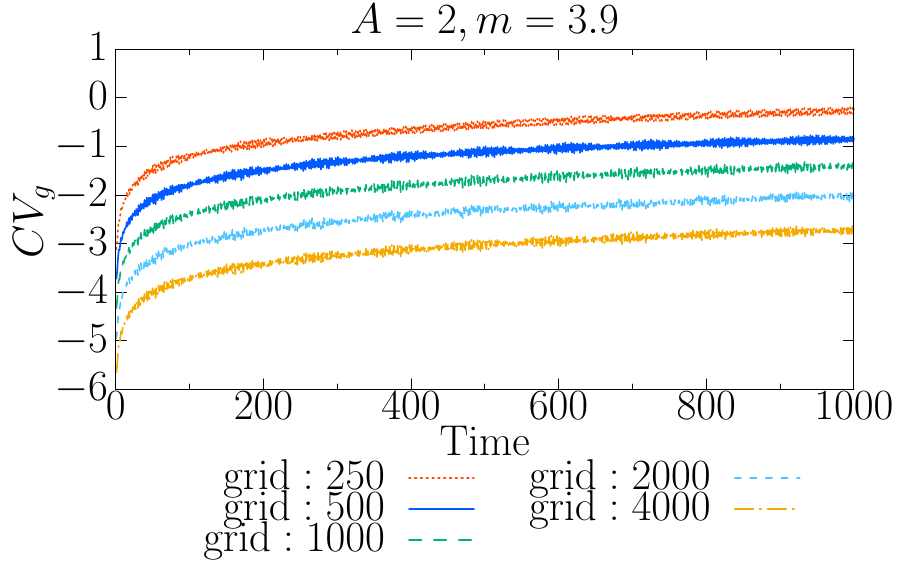}
  \end{minipage}
  \begin{minipage}{0.32\hsize}
    \includegraphics[width=\hsize]{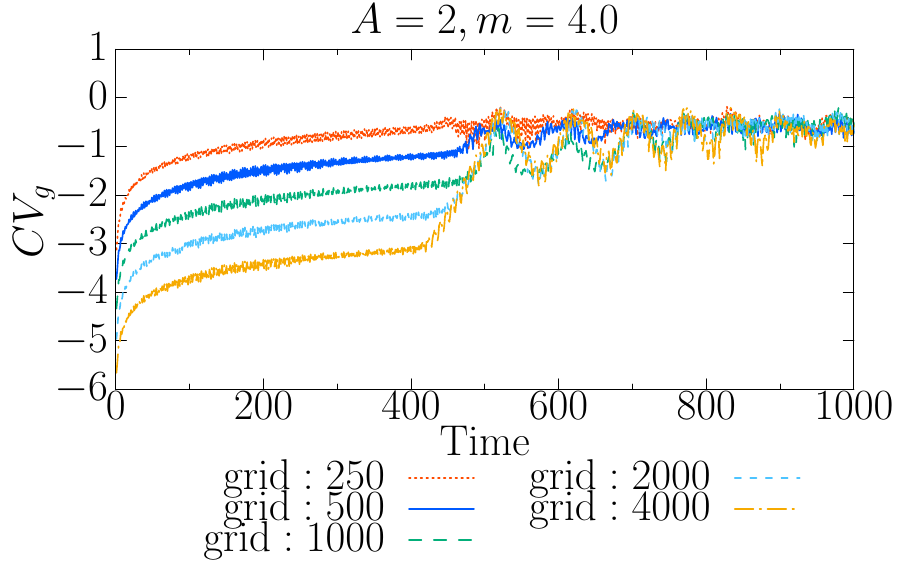}
  \end{minipage}
  \begin{minipage}{0.32\hsize}
    \includegraphics[width=\hsize]{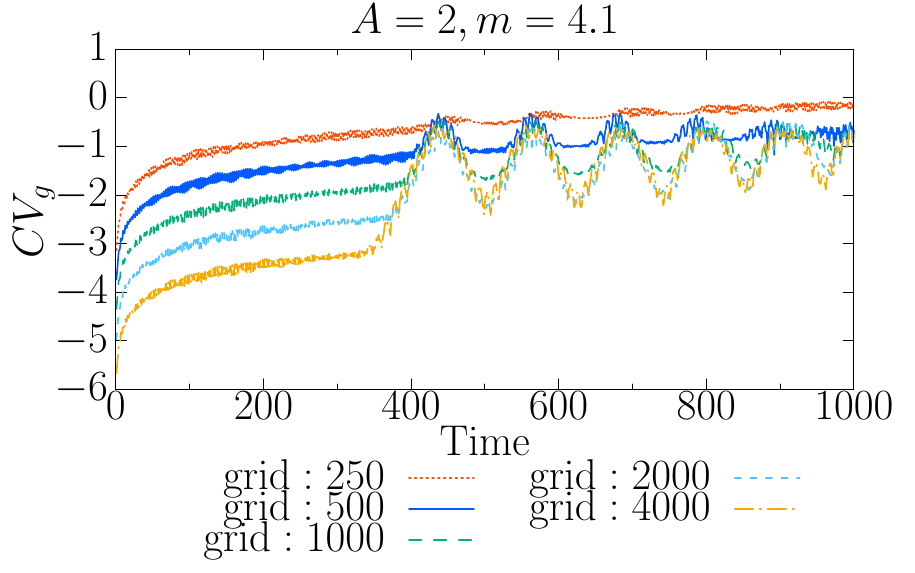}
  \end{minipage}
  \begin{minipage}{0.32\hsize}
    \includegraphics[width=\hsize]{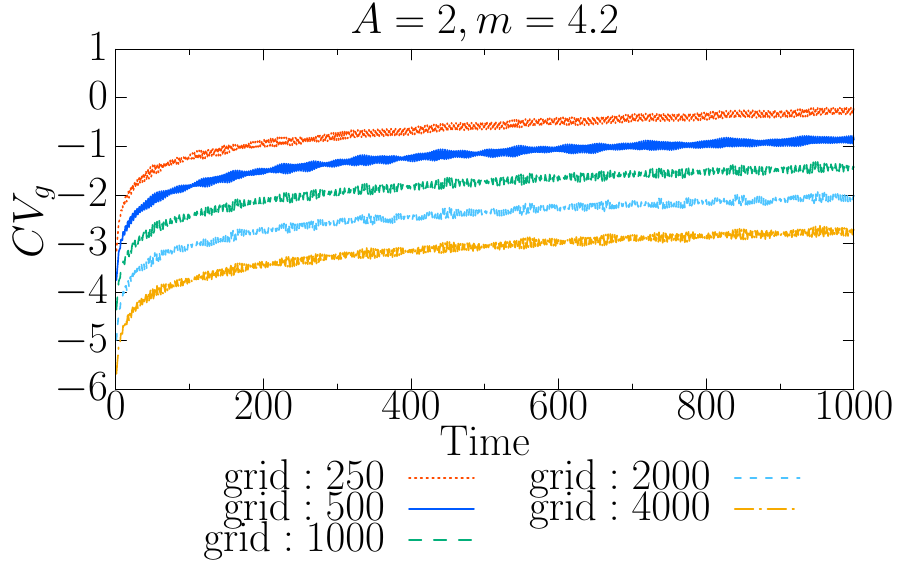}
  \end{minipage}
  \caption{
    Relative errors between $\phi$ with $8000$ grids and $\phi$ with other grid
    numbers when $A=2$ and $m=3.9$ to $4.2$.
    The vertical axis is $CV_g$ and the horizontal axis is time.
    The top-left panel is for $m=3.9$, the top-center one is for $m=4.0$, the
    top-right one is for $m=4.1$, and the bottom one is for $m=4.2$.
    The convergence seems not satisfied at either $t\geq 400$ for $m=4.0$ or
    $t\geq 350$ for $m=4.1$.
    \label{fig:H0ConvergenceA2_m39_42}
  }
\end{figure}
Fig. \ref{fig:H0ConvergenceA2_m39_42} shows the convergence of $\phi$ for $A=2$
and $m=3.9$ to $4.2$.
The convergence seems not satisfied for either $m=4$ or $m=4.1$.
\begin{figure}[htbp]
  \centering
  \begin{minipage}{0.32\hsize}
    \includegraphics[width=\hsize]{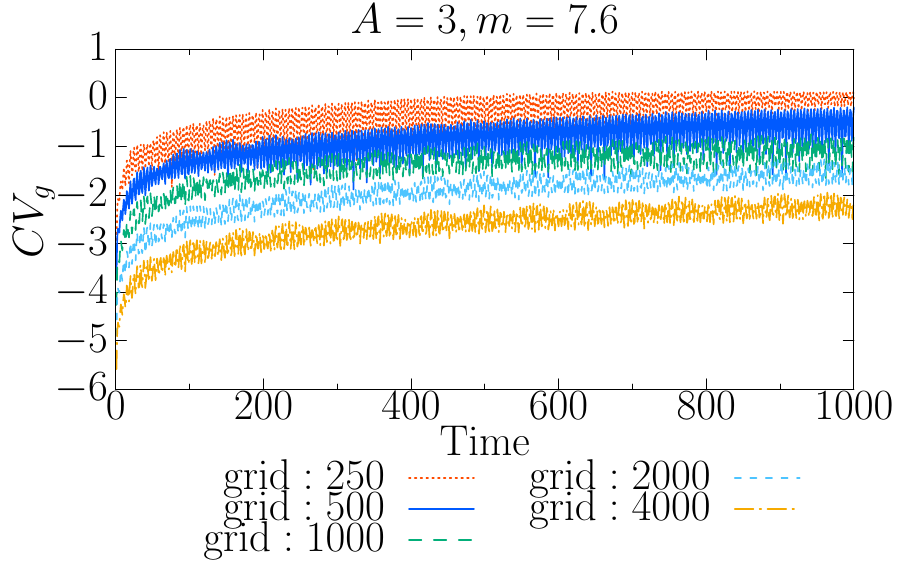}
  \end{minipage}
  \begin{minipage}{0.32\hsize}
    \includegraphics[width=\hsize]{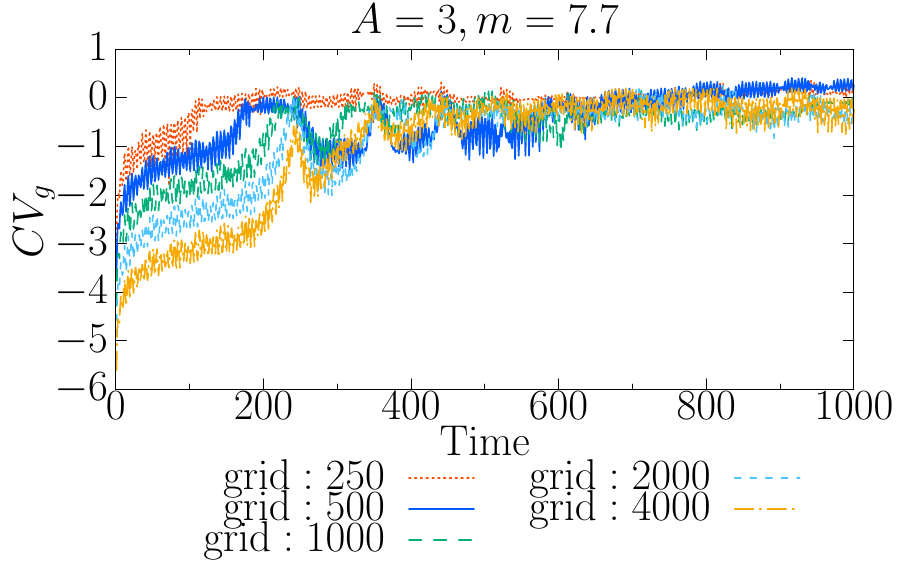}
  \end{minipage}
  \begin{minipage}{0.32\hsize}
    \includegraphics[width=\hsize]{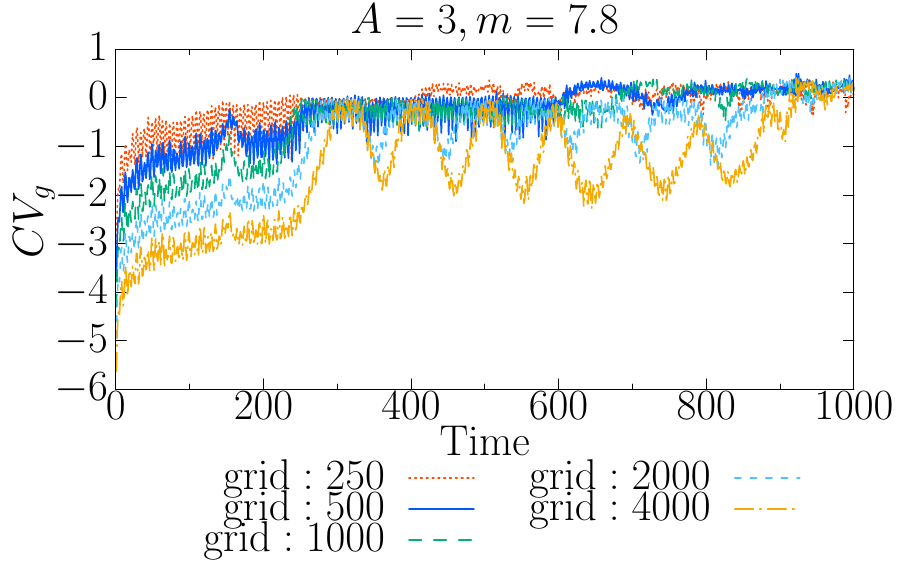}
  \end{minipage}
  \begin{minipage}{0.32\hsize}
    \includegraphics[width=\hsize]{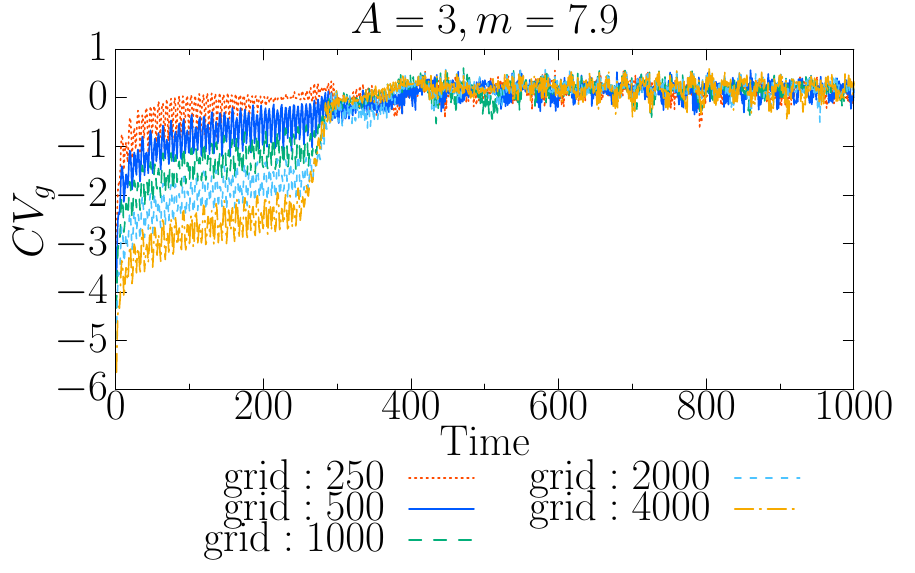}
  \end{minipage}
  \begin{minipage}{0.32\hsize}
    \includegraphics[width=\hsize]{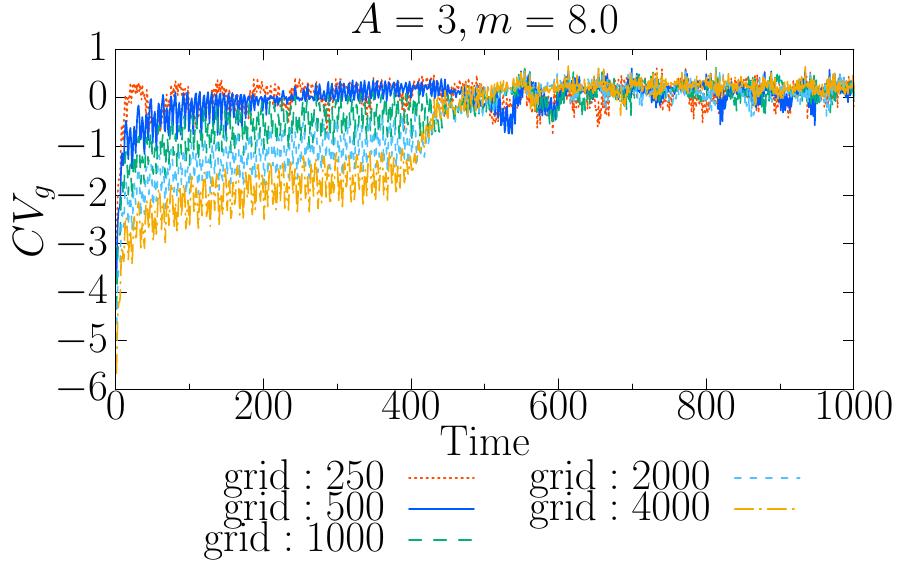}
  \end{minipage}
  \begin{minipage}{0.32\hsize}
    \includegraphics[width=\hsize]{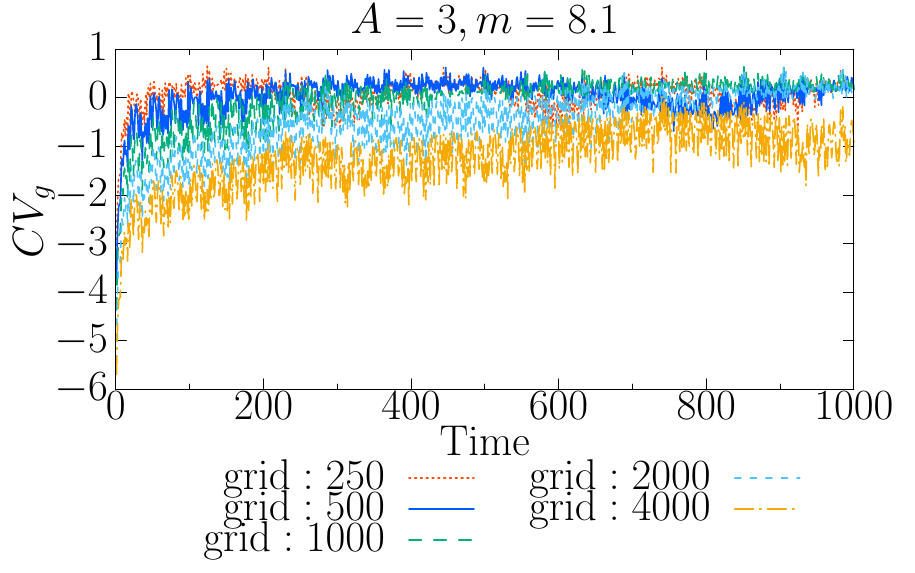}
  \end{minipage}
  \begin{minipage}{0.32\hsize}
    \includegraphics[width=\hsize]{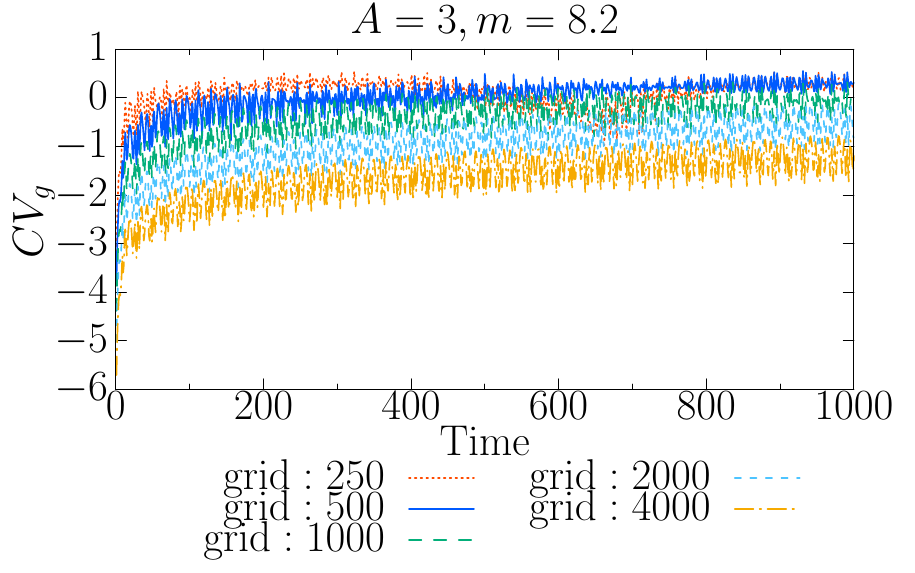}
  \end{minipage}
  \caption{
    Relative errors between $\phi$ with $8000$ grids and $\phi$ with other grid
    numbers when $A=3$ and $m=7.6$ to $8.2$.
    The top-left panel is for $m=7.6$, the top-center one is for $m=7.7$, the
    top-right one is for $m=7.8$, the center-left one is for $m=7.9$, the
    center one is for $m=8.0$, the center-right one is for $m=8.1$, and the
    bottom one is for $m=8.2$.
    The convergence seems not satisfied at $t\geq 200$ for $m=7.7$, at
    $t\geq 250$ for $m=7.8$, at $t\geq 250$ for $m=7.9$, at $t\geq 400$ for
    $m=8.0$, or at $t\geq 250$ for $m=8.1$.
    \label{fig:H0ConvergenceA3_m76_81}
  }
\end{figure}
On the other hand, Fig. \ref{fig:H0ConvergenceA3_m76_81} shows the convergence
for $A=3$ and $m=7.6$ to $8.2$.
The convergence seems not satisfied for $m=7.7$ to $8.1$.

To quantitatively evaluate convergence, we calculate $DCV_g(t)$ using \eqref{%
  eq:QuantitativeConv} at various $\varepsilon_{\mathrm{c}}$ values from $0.1$
to $0.4$.
Then, $\phi_{\mathrm{G}}(x)=\phi_{8000}(x)$ and $CV_{\bar{\mathrm{G}}}(t)
= CV_{4000}(t)$ since the maximum grid number is $8000$ and the second largest
grid number is $4000$ in these simulations.
We summarize the results of the convergence in Table \ref{tbl:convA2A3}.
\begin{table}[htbp]
  \centering
  \begin{tabular}{|c||c|c|c|c|}
    \hline
    \multicolumn{5}{|c|}{$A=2$}\\
    \hline
    \diagbox[width=1.1cm]{$\varepsilon_{\mathrm{c}}$}{$m$}
    & $3.9$ & $4.0$ & $4.1$ & $4.2$ \\
    \hline
    $0.1$ & $333$ & $423$ & $353$ & $450$\\
    \hline
    $0.15$ & $\bigcirc$ & $431$ & $354$ & $\bigcirc$\\
    \hline
    $0.2$ & $\bigcirc$ & $431$ & $354$ & $\bigcirc$\\
    \hline
    $0.25$ & $\bigcirc$ & $431$ & $363$ & $\bigcirc$\\
    \hline
    $0.3$ & $\bigcirc$ & $439$ & $363$ & $\bigcirc$\\
    \hline
    $0.35$ & $\bigcirc$ & $439$ & $365$ & $\bigcirc$\\
    \hline
    $0.4$ & $\bigcirc$ & $440$ & $365$ & $\bigcirc$\\
    \hline
  \end{tabular}
  \begin{tabular}{|c||c|c|c|c|c|c|c|}
    \hline
    \multicolumn{8}{|c|}{$A=3$}\\
    \hline
    \diagbox[width=1.1cm]{$\varepsilon_{\mathrm{c}}$}{$m$}
    & $7.6$ & $7.7$ & $7.8$ & $7.9$ & $8.0$ & $8.1$
    & $8.2$\\
    \hline
    $0.1$ & $68$ & $41$ & $42$ & $28$ & $11$ & $17$ & $37$\\
    \hline
    $0.15$ & $942$ & $212$ & $149$ & $255$ & $107$ & $109$ & $276$\\
    \hline
    $0.2$ & $\bigcirc$ & $222$ & $239$ & $255$ & $241$ & $164$ & $530$\\
    \hline
    $0.25$ & $\bigcirc$ & $227$ & $244$ & $255$ & $330$ & $164$ & $745$\\
    \hline
    $0.3$ & $\bigcirc$ & $231$ & $251$ & $260$ & $396$ & $219$ & $\bigcirc$\\
    \hline
    $0.35$ & $\bigcirc$ & $270$ & $255$ & $261$ & $396$ & $219$ & $\bigcirc$\\
    \hline
    $0.4$ & $\bigcirc$ & $273$ & $256$ & $261$ & $396$ & $267$ & $\bigcirc$\\
    \hline
  \end{tabular}
  \caption{
    Time when $DCV_{2000}>\varepsilon_{\mathrm{c}}$ for $A=2$ and $m=3.9$ to
    $4.2$, and for $A=3$ and $m=7.6$ to $8.2$.
    The left table is for $A=2$ and the right one is for $A=3$.
    $\bigcirc$ means that $DCV_{2000} \leq \varepsilon_{\mathrm{c}}$ is always
    satisfied at $0\leq t\leq 1000$.
    The values represent the times when $DCV_{2000}>\varepsilon_{\mathrm{c}}$.
    \label{tbl:convA2A3}
  }
\end{table}
$\bigcirc$ means that $DCV_{2000}\leq\varepsilon_{\mathrm{c}}$ is always
satisfied at $0\leq t\leq 1000$.
On the other hand, the values represent the times when $DCV_{2000} >
\varepsilon_{\mathrm{c}}$.

\section{Conclusion and discussion}

We showed some simulations of the semilinear Klein--Gordon equation with the
power-law nonlinear term in the flat spacetime.
The simulations were performed using the discrete equation, which was
constructed by a structure-preserving scheme for various mass $m$ values
from $3.9$ to $4.2$, where the amplitude of the initial value was $A=2$ and
$m$ ranged from $7.6$ to $8.2$ when $A=3$.
We proposed quantitative evaluation methods for stability and convergence.

For the threshold of stability, $\varepsilon_{\mathrm{s}}$, an appropriate value
must be selected so that the time when vibration occurs in Figs.
\ref{fig:WavesPhi8000_A2_m39_42}-\ref{fig:WavesPhi8000_A3_m76_82} and the time
when $SV_{8000}>\varepsilon_{\mathrm{s}}$ are the same, for each amplitude $A$.
In the case of $A=2$, from Fig. \ref{fig:WavesPhi8000_A2_m39_42},
the vibration appears to occur at $t \geq 500$ for $m=4.0$ and at $t \geq 700$
for $m=4.1$.
Therefore, $\varepsilon_{\mathrm{s}}$ is adopted $0.24$ for $A=2$.
In the case of $A=3$, from Fig. \ref{fig:WavesPhi8000_A3_m76_82}, the vibration
appears to occur at $t \geq 900$ for $m=7.8$, at $t \geq 300$ for $m=7.9$, and
at $t \geq 400$ for $m=8.0$.
Therefore, $\varepsilon_{\mathrm{s}}$ is also adopted $0.24$ for $A=3$.
On the other hand, for the threshold of convergence, $\varepsilon_{\mathrm{c}}$,
there are differences from $0.1$ to $0.4$ in Table \ref{tbl:convA2A3}.
$\varepsilon_{\mathrm{c}}$ is a value for quantitatively judging the
second-order convergence of the solutions.
For $A=2$ and $3$, there is no significant change in the number of times at
which $\varepsilon_{\mathrm{c}}\geq 0.15$ and $\varepsilon_{\mathrm{c}}\geq0.3$,
respectively.
Therefore, we adopt $0.15$ and $0.3$ as the thresholds for $A=2$ and $3$,
respectively.
Regarding $\varepsilon_{\mathrm{c}}$ for $A=3$ being greater than that for
$A=2$, this means that the convergence becomes worse as the amplitude of the
initial value increases.
It seems that the large initial amplitude due to the effect of nonlinear terms
worsens the convergence of the solutions.
Note that the appropriate values of $\varepsilon_{\mathrm{s}}$ and
$\varepsilon_{\mathrm{c}}$ depend on the parameters of the numerical
calculation, such as the amplitude of the initial value and the mass.
Thus, we have to investigate the appropriate values of thresholds under
different numerical calculation conditions.

In this study, we only investigated in the flat spacetime.
What we would like to investigate next in our future work is in a curved
spacetime.
The range of values for $m$ and $A$ is limited in this paper.
Future work will involve investigating the stability and convergence over a
wider range of values for $m$ and $A$.

\section*{Acknowledgments}
T.T. and M.N. were partially supported by JSPS KAKENHI Grant Numbers 21K03354
and 24K06855.
T.T. was partially supported by JSPS KAKENHI Grant Number 24K06856.


\end{document}